\documentclass[a4paper,11pt]{article}
\usepackage[T1]{fontenc}
\usepackage{lmodern}
\usepackage{amsmath,amssymb,amsthm,mathtools,mathrsfs}
\usepackage{microtype}
\usepackage{indentfirst} 
\usepackage{geometry}
\usepackage{needspace}
\usepackage[backend=biber,style=numeric,sorting=nyt,giveninits=true,doi=true,isbn=false,url=true]{biblatex}
\usepackage[hidelinks]{hyperref}
\numberwithin{equation}{section}
\newtheorem{theorem}{Theorem}[section]
\newtheorem{lemma}{Lemma}[section]
\newtheorem{corollary}{Corollary}[section]

\theoremstyle{definition}
\newtheorem{example}{Example}[section]
\theoremstyle{remark}
\newtheorem{remark}{Remark}[section]
\newcommand{\R}{\mathbb R}
\newcommand{\dd}{\,\mathrm d}
\DeclareMathOperator{\sgn}{sgn}
\hypersetup{pdfauthor={Ao Xiao},pdftitle={Multiple Positive Solutions of Minkowski-Curvature Neumann Problems with a Two-Step Indefinite Weight},pdfsubject={Constructive criteria for multiple positive solutions and a fixed-weight example with infinitely many solutions}}

\title{Multiple Positive Solutions of Minkowski-Curvature\\
Neumann Problems with a Two-Step Indefinite Weight}
\author{Ao Xiao\\[-1pt]
{\small Department of Mathematics, Northwest Normal University,}\\[-1pt]
{\small Lanzhou 730070, P.~R.~China}\\[-1pt]
{\small \href{mailto:xiaoao2math@163.com}{xiaoao2math@163.com}}}
\date{}
\begin{document}
\maketitle
\begin{abstract}
We study the multiplicity of positive solutions of the following
Minkowski-curvature Neumann problem
\[
\begin{aligned}
\left(\frac{u'}{\sqrt{1-(u')^2}}\right)'
+a(t)g(u)&=0,\quad t\in(0,T),\\
u'(0)=u'(T)&=0,
\end{aligned}
\]
where $a$ is a two-step sign-changing weight with negative mean.
We establish a multiplicity criterion for positive solutions, give a
negative answer to the uniqueness question raised by Boscaggin,
Feltrin and Zanolin, and prove the existence of infinitely many
positive solutions for a fixed weight. The proofs are based on the
shooting method and uniform asymptotic estimates.
\end{abstract}



\section{Introduction}\label{sec:introduction}

We consider the one-dimensional Minkowski-curvature Neumann problem
\begin{equation}\label{eq:problem}
 \begin{cases}
 \displaystyle\biggl(\frac{u'(t)}{\sqrt{1-(u'(t))^2}}\biggr)'
               +a(t)g(u(t))=0, & 0<t<T,\\[3pt]
 u'(0)=u'(T)=0,
 \end{cases}
\end{equation}
with the two-step weight
\begin{equation}\label{eq:weight}
 a(t)=\begin{cases}
 \alpha,&0<t<\tau,\\
 -\beta,&\tau<t<T=\tau+\sigma,
 \end{cases}
 \qquad \alpha,\beta,\tau,\sigma>0.
\end{equation}
Throughout the paper we assume
\begin{equation}\label{eq:g-assumption}
 g\in C^1([0,\infty)),\qquad
 g(s)>0,\quad g'(s)>0\quad(s>0).
\end{equation}
Set
\[
 \phi(v)=\frac{v}{\sqrt{1-v^2}},\qquad |v|<1.
\]
By a positive solution of \eqref{eq:problem} we mean a function
$u\in C^1([0,T])$ with $u>0$ and $|u'|<1$ on $[0,T]$ such that
$\phi(u')$ is absolutely continuous, the equation holds almost
everywhere and the boundary conditions are satisfied.
By \cite[Remark~2.1]{BFZ2021}, a necessary condition for the existence
of positive solutions is
\begin{equation}\label{eq:negative-mean}
 \int_0^T a(t)\dd t=\alpha\tau-\beta\sigma<0.
\end{equation}
Throughout, $\kappa=\sigma/\tau$ denotes the prescribed ratio of the
negative-step length to the positive-step length.

Boscaggin, Feltrin and Zanolin proved uniqueness for several step-weight
$\phi$-Laplacian problems, including power nonlinearities with
$p$-Laplacian fluxes \cite[Theorems~1.1--1.2]{BFZ2021}.
For the Minkowski-curvature problem with
$g_p(u)=e^{u^p}-1$, $p>1$, they established the existence of positive
solutions under Neumann and periodic boundary conditions
\cite[Theorem~2.1 and Corollary~2.1]{BFZ2023}. Their numerical results
suggested uniqueness \cite[Section~4.2, Eq.~(4.5)]{BFZ2021}.
The uniqueness question raised by Boscaggin, Feltrin and Zanolin is
recorded in the DEG1 Challenge \cite{DEG1Challenge}. It asks whether
a step weight with a single sign change and negative mean ensures
uniqueness of the positive solution for either boundary condition.
The special Neumann case $p=2$ with equal step lengths was studied
in \cite{Xiao2026}, where a counterexample and the exact multiplicity
pattern $1$--$2$--$3$--$2$--$1$ as the negative weight amplitude varies
were obtained for suitable parameters.
In this paper, for every $p>1$ and every $\kappa>0$, we construct a
weight for which the Neumann problem has at least three positive solutions
(Corollary~\ref{cor:superexp}).
Reflection gives a periodic weight with the same step-length ratio
and the same lower bound (Corollary~\ref{cor:periodic}). Thus we give
a negative answer to the uniqueness question for both boundary
conditions.

For semilinear indefinite Neumann problems, Feltrin and Sovrano
obtained three solutions by shooting \cite[Theorem~1.1]{FS2018},
and Tellini proved a high-multiplicity result for a superlinear
problem \cite{Tellini2018}. Further results on parameter-dependent
multiplicity can be found in \cite{FSHigh2018,FST2022}.
For the one-dimensional Minkowski-curvature Neumann equation,
Ma, Su and Zhao obtained at least two positive solutions for a general
indefinite weight and large coupling under their growth hypotheses
\cite{MSZ2023}, while Su, Ma and Yang studied $S$-shaped components
of positive Neumann solutions \cite{SMY2025}.
For radial problems, see \cite{BF2020,BCN2020}; for periodic and
subharmonic solutions, see \cite{BFPeriodic2020}; and for topological
methods in indefinite problems, see \cite{FeltrinBook2018}.

For one-dimensional Minkowski-curvature Dirichlet problems, Coelho,
Corsato, Obersnel and Omari obtained infinitely many positive solutions
under conditions near zero \cite{CCOO2012}. Lee, Sim and Yang studied
global bifurcation curves with repeated $S$-shapes and $\Sigma$-shapes
\cite{LSY2024}. In contrast to these approaches, we use the initial
height as the shooting parameter for a two-step Neumann weight.

The multiplicity criterion is expressed in terms of the profile
\[
 \ell(s)=\frac{g'(s)}{g(s)},\qquad
 F_\Lambda(s)=\ell(s)d(\Lambda g(s)),\qquad \Lambda=\alpha\tau,
\]
where $d$ is defined in \eqref{eq:d-ell-F}. In this paper, for each
prescribed $\kappa>0$, we take $T=(1+\kappa)L$ and both weight
amplitudes of order $L^{-1}$ as $L\downarrow0$; see
\eqref{eq:parameter-family}. The first-order residual formula holds
for every $g$ satisfying \eqref{eq:g-assumption} and every $\kappa>0$
(Lemma~\ref{lem:residual}). Strict sign alternation
of $F_\Lambda-\gamma$ for some $\gamma>0$ yields a lower bound on the
number of positive solutions which persists on an open set of
four-parameter weights (Theorem~\ref{thm:main}). We also obtain
criteria in terms of $\ell$ (Corollary~\ref{cor:simple}). In particular,
if $g'(0)=0$ and
\[
 \ell(\xi_2)<\min\{\ell(\xi_1),\ell(\xi_3)\}
 \qquad\text{for some }0<\xi_1<\xi_2<\xi_3,
\]
then suitable weights yield at least three positive solutions
(Corollary~\ref{cor:three}). These criteria allow $g(0)>0$.

Finally, for $g(s)=s^2e^{2s-\cos s}$, by establishing a uniform residual
estimate on a half-line
of initial heights (Lemma~\ref{lem:tail}), we construct, for each
prescribed ratio, a fixed weight admitting infinitely many positive
solutions tending uniformly to $+\infty$
(Theorem~\ref{prop:oscillatory}).

\section{Preliminaries}\label{sec:prelim}

The inverse Minkowski flux is $\psi(q)=q/\sqrt{1+q^2}$. We use
$|\psi|<1$ and $\psi'(q)=(1+q^2)^{-3/2}$. For $k\ge0$, $s>0$
and $\Lambda>0$, set
\begin{equation}\label{eq:d-ell-F}
 d(k)=2\int_0^1\frac{kr^2}{\sqrt{1+k^2r^2}}\dd r,
 \qquad \ell(s)=\frac{g'(s)}{g(s)},\qquad
 F_\Lambda(s)=\ell(s)d(\Lambda g(s)).
\end{equation}
The kernel $d$ and, in the case $g=g_2$ with equal step lengths,
the function $F_\Lambda$ agree with $W$ and $M_L$, respectively,
in \cite[Eqs.~(3.2)--(3.3)]{Xiao2026}, where the parameter $L$
denotes the integrated positive weight (denoted by $\Lambda$ here).

\begin{lemma}\label{lem:kernel}
The kernel $d$ is continuous and strictly increasing on $[0,\infty)$,
with $d(0)=0$, $0<d(k)<1$ and $d(k)\le 2k/3$ for $k>0$. Moreover,
\[
 d(k)=\tfrac23 k-\tfrac15 k^3+O(k^5)\quad(k\downarrow0),
 \qquad \lim_{k\to\infty}d(k)=1.
\]
\end{lemma}

\begin{proof}
By \eqref{eq:d-ell-F}, $d(k)=2\int_0^1r\psi(kr)\dd r$.
For $r>0$ the integrand $2r\psi(kr)$ increases strictly with $k$, and for $k>0$
it is bounded by both $2r$ and $2kr^2$, with strict inequality in
the first bound. Dominated convergence gives continuity and the
limit $2\int_0^1r\dd r=1$ at infinity. At zero, use
$(1+k^2r^2)^{-1/2}=1-k^2r^2/2+O(k^4)$ uniformly for $0\le r\le1$.
\end{proof}

Fix $\kappa,\Lambda,\gamma>0$ and put $S=1+\kappa$, $c_\kappa=S/2$.
For $L>0$, set
\begin{equation}\label{eq:parameter-family}
 \tau=L,\qquad\sigma=\kappa L,\qquad
 \alpha=\frac{\Lambda}{L},\qquad
 \beta=\frac{\Lambda(1+c_\kappa\gamma L)}{\kappa L}.
\end{equation}
Then $T=LS$, $\alpha\tau=\Lambda$ and
\begin{equation}\label{eq:mean-computation}
 \int_0^T a(t)\dd t=\alpha\tau-\beta\sigma
 =-\Lambda c_\kappa\gamma L<0.
\end{equation}
For an initial height $s>0$, set $x=t/L$, $U(x)=u(Lx)$ and
$q(x)=-\phi(u'(Lx))$. The identities $U_x=-L\psi(q)$ and
$q_x=La(Lx)g(U)$ give the shooting system
\begin{equation}\label{eq:scaled-IVP}
 \begin{cases}
 U_x=-L\psi(q),&0<x<S,\\[2pt]
 q_x=\Lambda g(U),&0<x<1,\\[2pt]
 q_x=-\dfrac{\Lambda(1+c_\kappa\gamma L)}{\kappa}g(U),
      &1<x<S,\\[3pt]
 U(0;s)=s,\quad q(0;s)=0.
 \end{cases}
\end{equation}
At $x=1$ the terminal data of the first step provide the initial data
for the second. The terminal
Neumann condition is equivalent to $q(S;s)=0$.

\begin{lemma}\label{lem:shooting-existence}
Let $I=[s_-,s_+]\subset(0,\infty)$, $s_-<s_+$. There exist
$L_I,Q_I>0$ such that, for $0<L<L_I$ and $s\in I$,
\eqref{eq:scaled-IVP} has a unique trajectory on $[0,S]$ with
\begin{equation}\label{eq:U-bounds}
 |U(x;s)-s|\le SL,
 \qquad U(x;s)\ge\frac{s_-}{2}>0\quad(0\le x\le S).
\end{equation}
Moreover, $|q(x;s)|\le Q_I$, and for each such $L$ the trajectory
depends continuously on $(x,s)\in[0,S]\times I$.
\end{lemma}

\begin{proof}
Take $L_I\le1$ with $SL_I<\min\{s_-/2,1\}$. The two vector fields
are locally Lipschitz on $(0,\infty)\times\R$. As long as the
trajectory exists, $|U(x;s)-s|\le L\int_0^x|\psi(q(y;s))|\dd y\le SL$.
Thus $U\in J=[s_-/2,s_++1]$. With $M_0=\max_J g$, integration
up to $x\le S$ gives
\[
 |q(x;s)|\le\Lambda M_0\min\{x,1\}
 +\frac{\Lambda(1+c_\kappa\gamma L)}{\kappa}M_0(x-1)_+
 \le \Lambda M_0(2+c_\kappa\gamma)=:Q_I,
\]
where $(x-1)_+=\max\{x-1,0\}$. The trajectory stays in the compact
rectangle $J\times[-Q_I,Q_I]$, so continuation on each step and
matching at $x=1$ give existence up to $S$. Uniqueness and continuous
dependence follow from local Lipschitz continuity on each step.
\end{proof}

Whenever a positive trajectory satisfies $q(S;s)=0$, the function
$u(t)=U(t/L;s)$ solves \eqref{eq:problem}. Indeed,
\begin{equation}\label{eq:uj-derivative}
 u'(t)=-\psi(q(t/L;s)),\qquad
 \frac{\mathrm d}{\mathrm dt}\phi(u'(t))
 =-L^{-1}q_x(t/L;s)=-a(t)g(u(t))\quad\text{a.e.}
\end{equation}
The flux is continuous and piecewise $C^1$ with bounded derivative,
hence absolutely continuous. The first identity gives $u\in C^1$,
$|u'|<1$ and $u'(0)=u'(T)=0$. Also, $q$ increases strictly from
zero on the first step and decreases strictly to zero on the second.
Thus $u'<0$ in $(0,T)$.

Define
\begin{equation}\label{eq:A-definition}
 A_\kappa(x)=\begin{cases}
 x,&0\le x\le1,\\[2pt]
 \dfrac{S-x}{\kappa},&1\le x\le S.
 \end{cases}
\end{equation}
Then $A_\kappa$ is absolutely continuous, vanishes at both endpoints,
and has derivative $1$ on $(0,1)$ and $-1/\kappa$ on $(1,S)$.
The normalized terminal residual is
\begin{equation}\label{eq:R-definition}
 R_{\Lambda,L,\kappa}(s)
 =\frac{q(S;s)}{\Lambda Lg(s)}.
\end{equation}
We suppress the fixed parameter $\gamma$ in this notation.
Since $\Lambda Lg(s)>0$, $R_{\Lambda,L,\kappa}(s)$ has the sign
of the terminal flux and its zeros give Neumann solutions by
\eqref{eq:uj-derivative}.
\begin{lemma}
\label{lem:residual}
Under the hypotheses of Lemma~\ref{lem:shooting-existence}, the
residual \eqref{eq:R-definition} satisfies, for $s\in I$ and
$0<L<L_I$,
\begin{align}\label{eq:exact-identity}
 R_{\Lambda,L,\kappa}(s)
 ={}&\int_0^S A_\kappa(x)
       \frac{g'(U(x;s))}{g(s)}\psi(q(x;s))\dd x
 \notag\\
 &-\frac{c_\kappa\gamma}{\kappa}
       \int_1^S\frac{g(U(x;s))}{g(s)}\dd x.
\end{align}
In addition,
\begin{equation}\label{eq:uniform-limit}
 \lim_{L\downarrow0}
 \sup_{s\in I}\bigl|
 R_{\Lambda,L,\kappa}(s)
 -c_\kappa(F_\Lambda(s)-\gamma)\bigr|=0.
\end{equation}
\end{lemma}

\begin{proof}
With the triangular function $A=A_\kappa$ from \eqref{eq:A-definition},
the flux equation reads
\begin{equation}\label{eq:q-compact}
 q_x(x;s)=\Lambda A_\kappa'(x)g(U(x;s))
  -\frac{\Lambda c_\kappa\gamma L}{\kappa}
       \mathbf1_{(1,S)}(x)g(U(x;s)).
\end{equation}
Both $A$ and $g\circ U$ are absolutely continuous. Since
$A(0)=A(S)=q(0;s)=0$ and $U_x=-L\psi(q)$, integration by parts gives
\begin{align*}
 R_{\Lambda,L,\kappa}(s)
 &=\frac{1}{Lg(s)}\int_0^S A'g(U)\dd x
    -\frac{c_\kappa\gamma}{\kappa}\int_1^S\frac{g(U)}{g(s)}\dd x\\
 &=-\frac{1}{Lg(s)}\int_0^S A g'(U)U_x\dd x
    -\frac{c_\kappa\gamma}{\kappa}\int_1^S\frac{g(U)}{g(s)}\dd x\\
 &=\int_0^S A\frac{g'(U)}{g(s)}\psi(q)\dd x
    -\frac{c_\kappa\gamma}{\kappa}\int_1^S\frac{g(U)}{g(s)}\dd x.
\end{align*}

For the limit, put $q_0(x;s)=\Lambda g(s)A(x)$,
$M_1=\max_J|g'|$ and $m=\min_I g>0$, with $J$ and $M_0$ as in the
proof of Lemma~\ref{lem:shooting-existence}. By \eqref{eq:U-bounds},
$|g(U)-g(s)|\le M_1SL$. Since $q(0;s)=q_0(0;s)=0$, integrating
\eqref{eq:q-compact} gives
\begin{align*}
 \|q-q_0\|_\infty
 &\le \Lambda M_1SL\int_0^S|A'(x)|\dd x
   +\frac{\Lambda c_\kappa\gamma LM_0}{\kappa}\int_1^S\dd x\\
 &=\Lambda(2M_1S+c_\kappa\gamma M_0)L=:C_q L,
\end{align*}
where the supremum norm is taken over $[0,S]\times I$.

Let $\omega_{g'}(\delta)=\sup\{|g'(v)-g'(w)|:v,w\in J,\
|v-w|\le\delta\}$. Since $g'$ is continuous on $J$,
$\omega_{g'}(\delta)\to0$ as $\delta\downarrow0$. Using
$|\psi|,|\psi'|\le1$, we obtain
\begin{align*}
 \left|\frac{g'(U)}{g(s)}\psi(q)-\ell(s)\psi(q_0)\right|
 &\le \frac{|g'(U)-g'(s)|+M_1|q-q_0|}{m}
 \le\frac{\omega_{g'}(SL)+M_1C_qL}{m},\\
 \left|\frac1\kappa\int_1^S\frac{g(U)}{g(s)}\dd x-1\right|
 &\le\frac1\kappa\int_1^S\frac{M_1SL}{m}\dd x
 =\frac{M_1SL}{m}.
\end{align*}
For $k>0$, substituting $r=(S-x)/\kappa$ on the second step gives
\begin{equation}\label{eq:triangle-integral}
\begin{aligned}
 \int_0^S A_\kappa(x)\psi(kA_\kappa(x))\dd x
 &=\int_0^1 x\psi(kx)\dd x
   +\int_1^S\frac{S-x}{\kappa}
       \psi\!\left(k\frac{S-x}{\kappa}\right)\dd x\\
 &=(1+\kappa)\int_0^1 r\psi(kr)\dd r\\
 &=c_\kappa d(k).
\end{aligned}
\end{equation}
Take $k=\Lambda g(s)$. In \eqref{eq:exact-identity} the first
integral differs from $c_\kappa F_\Lambda(s)$ by at most
$S[\omega_{g'}(SL)+M_1 C_qL]/m$, since $0\le A\le1$. Together
with the second integral estimate this yields
\[
 \sup_{s\in I}|R_{\Lambda,L,\kappa}(s)-c_\kappa(F_\Lambda(s)-\gamma)|
 \le \frac{S}{m}\bigl[\omega_{g'}(SL)+M_1C_qL\bigr]
       +\frac{c_\kappa\gamma M_1SL}{m}\longrightarrow0 .\qedhere
\]
\end{proof}

\section{Existence of multiple positive solutions}\label{sec:proof}

\begin{theorem}\label{thm:main}
Assume \eqref{eq:g-assumption} and fix an integer $N\ge1$ and a ratio
$\kappa>0$. Suppose there are $\Lambda,\gamma>0$ and heights
$0<s_0<\cdots<s_N$ such that
\begin{equation}\label{eq:main-alternation}
 \bigl(F_\Lambda(s_j)-\gamma\bigr)
 \bigl(F_\Lambda(s_{j+1})-\gamma\bigr)<0,
 \qquad j=0,\ldots,N-1.
\end{equation}
Then for every sufficiently small $L>0$, the weight
\eqref{eq:parameter-family} has negative mean and \eqref{eq:problem}
has at least $N$ distinct positive solutions. For each such $L$, both
conclusions persist in an open neighborhood of $(\alpha,\beta,\tau,\sigma)$
in $(0,\infty)^4$.
\end{theorem}

\begin{proof}
Set $I=[s_0,s_N]$ and $\Delta_j=F_\Lambda(s_j)-\gamma$.
The weight has negative integral by \eqref{eq:mean-computation}.
Lemma~\ref{lem:shooting-existence} gives positive trajectories and
continuity of the terminal residual on $I$. By \eqref{eq:uniform-limit}
we may choose $0<L_0\le L_I$ so that for $0<L<L_0$,
\[
 \sup_{s\in I}|R_{\Lambda,L,\kappa}(s)-c_\kappa(F_\Lambda(s)-\gamma)|
 <\frac{c_\kappa}{2}\min_{0\le j\le N}|\Delta_j|.
\]
Hence $\sgn R_{\Lambda,L,\kappa}(s_j)=\sgn\Delta_j$, and consecutive
residuals have opposite signs by \eqref{eq:main-alternation}. The
intermediate value theorem gives a zero $\widehat s_j\in(s_j,s_{j+1})$
for $0\le j<N$. The corresponding trajectories give positive Neumann
solutions by \eqref{eq:uj-derivative}, distinct because their initial
heights lie in disjoint intervals.

Fix such an $L$ and denote the resulting parameter vector by $\theta_0$.
For $\theta=(\alpha,\beta,\tau,\sigma)$ near $\theta_0$, keep the
interface fixed by taking $t=\tau y$ on $[0,1]$ and
$t=\tau+\sigma(y-1)$ on $[1,2]$. The height and negative flux in these
coordinates satisfy
\begin{equation}\label{eq:parameter-shooting}
 \begin{aligned}
 (Y_y,Q_y)&=
 \begin{cases}
 (-\tau\psi(Q),\,\alpha\tau g(Y)),&0<y<1,\\
 (-\sigma\psi(Q),\,-\beta\sigma g(Y)),&1<y<2,
 \end{cases}\\
 (Y(0),Q(0))&=(s,0).
 \end{aligned}
\end{equation}
In particular, $Q(2;s,\theta)$ is the original terminal negative flux.
At $\theta_0$ it has the sign of $R_{\Lambda,L,\kappa}(s)$, and the
base trajectories from $I$ lie in
$K_0=[s_0/2,s_N+1]\times[-Q_I,Q_I]$.
Enlarge this rectangle to
\[
 K=[s_0/4,s_N+2]\times[-Q_I-1,Q_I+1]\Subset(0,\infty)\times\R .
\]
On $K$ and a compact parameter neighborhood of $\theta_0$, both
fields in \eqref{eq:parameter-shooting} have bounded first derivatives
in the states and parameters. For trajectories with the same initial
height, set
\[
 D(y)=|Y_\theta(y)-Y_{\theta_0}(y)|
       +|Q_\theta(y)-Q_{\theta_0}(y)|.
\]
Until the first exit from $K$, integration across the two steps gives
\[
 D(y)\le C_K y|\theta-\theta_0|+C_K\int_0^y D(z)\dd z.
\]
By Gronwall's inequality,
\[
 D(y)\le 2C_K e^{2C_K}|\theta-\theta_0|,\qquad 0\le y\le2,
\]
as long as the perturbed trajectory remains in $K$, with $C_K$
independent of $s\in I$. Choose $\theta$ sufficiently close to
$\theta_0$ so that
\[
 2C_K e^{2C_K}|\theta-\theta_0|
 <\tfrac12\min\{s_0/4,1\}.
\]
Since $\operatorname{dist}(K_0,\partial K)\ge\min\{s_0/4,1\}$,
a first exit is impossible. Continuation on each step gives existence
and positivity on $[0,2]$ for all $s\in I$. Continuous dependence on
the initial data on each step gives continuity in $s$ as well.

At the finitely many test heights the base terminal fluxes are nonzero.
By continuity there is an open neighborhood $\mathcal U$ of $\theta_0$
on which
\[
 Q(2;s_j,\theta)Q(2;s_{j+1},\theta)<0\quad(0\le j<N),
 \qquad \alpha\tau-\beta\sigma<0.
\]
For every $\theta\in\mathcal U$, each interval $(s_j,s_{j+1})$ again
contains a zero. Inverting the two changes of variable in
\eqref{eq:parameter-shooting} gives $N$ distinct positive Neumann
solutions. This argument requires no simplicity of the zeros.
\end{proof}

\begin{corollary}\label{cor:simple}
Assume \eqref{eq:g-assumption}.  Let $N\ge1$ be an integer and $\kappa>0$.
If there exist $\gamma>0$ and $0<s_0<\cdots<s_N$ such that
\begin{equation}\label{eq:ell-alternation}
 \bigl(\ell(s_j)-\gamma\bigr)
 \bigl(\ell(s_{j+1})-\gamma\bigr)<0,
 \qquad j=0,\ldots,N-1,
\end{equation}
then a strictly negative-mean weight \eqref{eq:weight} with
$\sigma/\tau=\kappa$ can be chosen so that \eqref{eq:problem}
has at least $N$ distinct positive solutions.
\end{corollary}

\begin{proof}
At the finitely many test heights $g(s_j)>0$, so
Lemma~\ref{lem:kernel} gives $F_\Lambda(s_j)\to\ell(s_j)$ as
$\Lambda\to\infty$. Choose one $\Lambda$ for which
\[
 \max_{0\le j\le N}|F_\Lambda(s_j)-\ell(s_j)|
 <\tfrac12\min_{0\le j\le N}|\ell(s_j)-\gamma|.
\]
Then \eqref{eq:ell-alternation} implies \eqref{eq:main-alternation}.
Theorem~\ref{thm:main} applies for all sufficiently small $L>0$.
\end{proof}

\begin{corollary}\label{cor:three}
Assume \eqref{eq:g-assumption}, $g'(0)=0$, and suppose that
there exist $0<\xi_1<\xi_2<\xi_3$ for which
$\ell(\xi_2)<\min\{\ell(\xi_1),\ell(\xi_3)\}$.
Then, for every prescribed $\kappa>0$, there is a weight
\eqref{eq:weight} with $\sigma/\tau=\kappa$ and
$\alpha\tau-\beta\sigma<0$ such that \eqref{eq:problem} has at
least three distinct positive solutions. The conclusion persists in an
open neighborhood in the four-parameter space around each weight
constructed for sufficiently small $L>0$.
\end{corollary}

\begin{proof}
Choose $\ell(\xi_2)<\gamma<\min\{\ell(\xi_1),\ell(\xi_3)\}$.
By Lemma~\ref{lem:kernel}, $F_\Lambda(\xi_i)\to\ell(\xi_i)>\gamma$
for $i=1,3$ as $\Lambda\to\infty$. Choose $\Lambda$ large enough for
both inequalities $F_\Lambda(\xi_i)>\gamma$ to hold. At $\xi_2$,
$d<1$ gives $F_\Lambda(\xi_2)<\ell(\xi_2)<\gamma$.
For this fixed $\Lambda$ the bound $d(k)\le 2k/3$ yields
\[
 0<F_\Lambda(s)=\frac{g'(s)}{g(s)}d(\Lambda g(s))
 \le\frac{2\Lambda}{3}g'(s)\longrightarrow0\qquad(s\downarrow0),
\]
since $g'(0)=0$. Choose $s_0\in(0,\xi_1)$ with
$F_\Lambda(s_0)<\gamma$. The signs at $s_0,\xi_1,\xi_2,\xi_3$ are
then $-,+,-,+$, and Theorem~\ref{thm:main} gives the conclusion.
\end{proof}

\begin{corollary}\label{cor:superexp}
For every $p>1$ and every prescribed $\kappa>0$, there is a weight
\eqref{eq:weight} of negative mean and length ratio $\kappa$ for which
\begin{equation}\label{eq:superexp-Neumann}
 \begin{cases}
 \displaystyle\biggl(\frac{u'}{\sqrt{1-(u')^2}}\biggr)'
        +a(t)(e^{u^p}-1)=0, &0<t<T,\\
 u'(0)=u'(T)=0,
 \end{cases}
\end{equation}
has at least three distinct positive solutions. The negative mean and
the lower bound of three persist in an open neighborhood of the
constructed four-parameter weight.
\end{corollary}

\begin{proof}
Set $g_p(u)=e^{u^p}-1$ for $u\ge0$.
For $u>0$, $g_p(u)>0$ and $g_p'(u)=pu^{p-1}e^{u^p}>0$.
Since $g_p(u)\sim u^p$ and $g_p'(u)\to0$ as $u\downarrow0$,
we have $g_p\in C^1([0,\infty))$ and $g_p(0)=g_p'(0)=0$.
For $u>0$ its logarithmic derivative is
\[
 \ell_p(u)=\frac{pu^{p-1}e^{u^p}}{e^{u^p}-1}
          =\frac{pu^{p-1}}{1-e^{-u^p}} .
\]
Set $h(z)=z/(e^z-1)$ for $z>0$. Logarithmic differentiation gives
\begin{align*}
 \frac{u\ell_p'(u)}{\ell_p(u)}
 &= (p-1)-\frac{pu^pe^{-u^p}}{1-e^{-u^p}}
  =(p-1)-p h(u^p),\\
 h'(z)&=\frac{e^z(1-z)-1}{(e^z-1)^2}<0 ,
\end{align*}
where the strict inequality follows from $1-z<e^{-z}$ for $z>0$.
Since $h(0^+)=1$ and $h(+\infty)=0$, the expression
$(p-1)-p h(u^p)$ increases strictly from $-1$ to $p-1>0$.
Consequently, there is a unique $b_p>0$ such that $\ell_p'(b_p)=0$,
with $\ell_p'<0$ on $(0,b_p)$ and $\ell_p'>0$ on $(b_p,\infty)$.
Thus $\ell_p$ attains its strict global minimum at $b_p$.
Any $\xi_1<b_p<\xi_3$ with $\xi_1>0$ therefore
satisfies the hypothesis of Corollary~\ref{cor:three} with $\xi_2=b_p$.
\end{proof}

\Needspace{8\baselineskip}
\begin{example}\label{ex:numerical}
Take $p=3/2$, $\kappa=2$, $\Lambda=30$, $\gamma=2.25$ and
$L=0.001$. Formula \eqref{eq:parameter-family} gives
\[
 (\alpha,\beta,\tau,\sigma)
   =(30000,\,15050.625,\,0.001,\,0.002).
\]
Here $S=3$ and $SL=0.003$. For $s\in[0.004,2.8]$, the speed
bound gives $U(x;s)\ge s-SL\ge0.001$ as long as the trajectory
exists. Since also $U\le s+SL$ and the flux is bounded by integrating
its equation, continuation and continuous dependence on both steps
hold as in Lemma~\ref{lem:shooting-existence}. In particular, the
residual is continuous on this interval, and the integration by parts
in Lemma~\ref{lem:residual} gives \eqref{eq:exact-identity} for this
fixed $L$.

To enclose the residual, put $B=1+c_\kappa\gamma L$ and fix a test
height $s$. If $|U-s|\le\delta\le SL$, then the monotonicity of $g$
and the two flux equations give $|q|\le\Lambda B g(s+\delta)$.
Thus $|U-s|\le\delta_5$, where
\[
 \delta_0=SL,\qquad
 \delta_{j+1}=SL\psi\bigl(\Lambda B g(s+\delta_j)\bigr),
 \quad j=0,\ldots,4.
\]
With $m=g(s-\delta_5)$ and $M=g(s+\delta_5)$, for $0\le r\le1$ we have
\begin{align*}
 \Lambda mr&\le q(r;s)\le\Lambda Mr,\\
 \Lambda(m-BMr)&\le q(1+\kappa r;s)\le\Lambda(M-Bmr).
\end{align*}
Both $g$ and $g'$ are increasing for $p=3/2$. Set
\[
 H=\frac{[g'(s-\delta_5),g'(s+\delta_5)]}{g(s)},\qquad
 G=\frac{[m,M]}{g(s)}.
\]
For $n=2000$, let $r_i=i/n$ and $J_i=[r_i,r_{i+1}]$.
Using $A_\kappa(r)=r$ and $A_\kappa(1+\kappa r)=1-r$ in
\eqref{eq:exact-identity}, interval arithmetic gives
\begin{align*}
 R_{\Lambda,L,\kappa}(s)\in{}&\frac1n\sum_{i=0}^{n-1}
 \bigl(J_iH\psi(Q_i^+)+\kappa(1-J_i)H\psi(Q_i^-)\bigr)
 -c_\kappa\gamma G,\\
 Q_i^+={}&[\Lambda mr_i,\Lambda Mr_{i+1}],\\
 Q_i^-={}&[\Lambda(m-BMr_{i+1}),\Lambda(M-Bmr_i)],
\end{align*}
where all operations are interval operations and $\psi([a,b])=
[\psi(a),\psi(b)]$. To evaluate the bounds, all interval endpoints
are rounded outwards to multiples of $10^{-30}$ (in particular,
$\delta_j$ is replaced by an upper bound at each iteration).
Square roots are enclosed by integer square roots. For $0\le z\le5$,
$e^z$ is enclosed using the Taylor sum of degree $80$ and a
remainder bounded above by $2z^{81}/81!$. All exponential arguments
in this computation lie in $[0,5]$. The resulting enclosures, rounded
outwards to three decimal places, are
\[
 \begin{array}{c@{\qquad}c}
 s & \text{Enclosure of }R_{\Lambda,L,\kappa}(s)\\ \hline
 0.004 & [-0.612,-0.450]\\
 0.08  & [7.355,8.943]\\
 1.6   & [-0.143,-0.050]\\
 2.8   & [0.365,0.484]
 \end{array}
\]
The enclosures give the strict sign pattern $-,+,-,+$ at the four
test heights. By continuity, the residual has a zero in each of
the three consecutive open intervals. Each zero gives a positive
solution of \eqref{eq:superexp-Neumann}, and the solutions are
distinct because their initial heights lie in disjoint intervals.

For comparison, integrating the two steps of \eqref{eq:scaled-IVP}
in succession with SciPy's DOP853 solver (relative tolerance
$10^{-12}$ and absolute tolerance $10^{-14}$ for both components)
gives the numerical residuals
\[
 \begin{array}{c@{\qquad}rrrr}
 s & 0.004 & 0.08 & 1.6 & 2.8\\ \hline
 R_{\Lambda,L,\kappa}(s) & -0.5214 & 8.1766 & -0.0910 & 0.4285
 \end{array}
\]
to four decimal places. Brent's method on the three consecutive
intervals gives the approximate initial values
$0.005594$, $1.197547$ and $1.970772$. These approximations are not
used in the existence proof.
\end{example}

\begin{corollary}\label{cor:periodic}
For every $p>1$ and every prescribed $\kappa>0$, there exists a
two-step weight with one sign change in the interior of the
period interval, strictly negative integral and ratio $\kappa$
of the negative-step length to the positive-step length for
which the equation in
\eqref{eq:superexp-Neumann}, with periodic boundary conditions,
has at least three distinct positive solutions.
\end{corollary}

\begin{proof}
Take three solutions $u_1,u_2,u_3$ for the same Neumann weight
$a$ on $[0,T]$ from Corollary~\ref{cor:superexp}. Following \cite[Section~2.2]{BFZ2021},
set $w_j(t)=u_j(|t|)$ and $\widetilde a(t)=a(|t|)$ on $[-T,T]$
and extend them $2T$-periodically. The zero derivatives at $0$
and $T$ ensure that $w_j$ is $C^1$ across the reflection points
and that $\phi(w_j')$ is absolutely continuous on each period.
Since $\phi$ is odd, the equation with $\widetilde a$ holds
almost everywhere on the reflected interval and hence on each
period.

On translating $[-\tau,2T-\tau]$ to $[0,2T]$, the weight is
$\alpha$ on $(0,2\tau)$ and $-\beta$ on $(2\tau,2T)$. Its
integral is $2(\alpha\tau-\beta\sigma)<0$, and its length ratio
is $\sigma/\tau=\kappa$. The translated solutions satisfy the
periodic conditions and remain distinct, since at $t=\tau$
their values are the distinct initial heights $u_j(0)$.
\end{proof}

Thus, Corollaries~\ref{cor:superexp} and~\ref{cor:periodic} give a
negative answer to the uniqueness question raised by Boscaggin,
Feltrin and Zanolin \cite{BFZ2021,DEG1Challenge}. For every $p>1$ and
prescribed step-length ratio $\kappa>0$, suitable two-step weights
with negative mean yield at least three positive solutions for both
Neumann and periodic boundary value problems.

\begin{remark}\label{rem:regime}
For fixed $\Lambda,\gamma,\kappa$, the weights in
\eqref{eq:parameter-family} have $T=(1+\kappa)L\to0$ and amplitudes
of order $L^{-1}$. As $L\downarrow0$,
\[
 \frac{\beta\sigma}{\alpha\tau}=1+c_\kappa\gamma L\to1^+,\qquad
 \int_0^T a(t)\dd t=-\Lambda c_\kappa\gamma L\to0^-,
 \qquad \frac1T\int_0^T a(t)\dd t=-\frac{\Lambda\gamma}{2}.
\]
Thus the two step integrals approach balance while the average remains
strictly negative. The height changes by at most $SL$ by
\eqref{eq:U-bounds}. These scalings describe the sufficient
construction used here.
\end{remark}

\section{Existence of infinitely many positive solutions}\label{sec:simple}

For a finite set of sign changes, Theorem~\ref{thm:main} chooses $L$
after the extreme test heights have been fixed. To use one weight for
an unbounded sequence of heights, we estimate the residual uniformly
on $[2,\infty)$. Consider
\begin{equation}\label{eq:example-g}
 g(s)=s^2e^{2s-\cos s},\qquad s\ge0.
\end{equation}
This function belongs to $C^1([0,\infty))$, and for $s>0$,
\[
 g'(s)=e^{2s-\cos s}\bigl[2s+s^2(2+\sin s)\bigr]>0,\qquad
 \ell(s)=2/s+2+\sin s .
\]
At $s=0$ we have $g'(0)=0$. Thus \eqref{eq:g-assumption} holds, and both $\ell$ and
$\ell'(s)=-2/s^2+\cos s$ are bounded on $[1,\infty)$.

\begin{lemma}\label{lem:tail}
Let $g$ be given by \eqref{eq:example-g}. Fix $\kappa,\Lambda>0$
and set $\gamma=2$. There exist $L_*>0$ and $C>0$ such that, for
every $0<L<L_*$, all trajectories of \eqref{eq:scaled-IVP} with
$s\ge2$ exist and stay positive on $[0,S]$. Their terminal residual
is continuous in $s$ and satisfies
\begin{equation}\label{eq:oscillatory-uniform-error}
 \sup_{s\ge2}\bigl|R_{\Lambda,L,\kappa}(s)-c_\kappa(F_\Lambda(s)-2)\bigr|
 \le CL .
\end{equation}
\end{lemma}

\begin{proof}
Choose $L_*\le1/S$. For $s\ge2$ the speed bound gives
$|U(x;s)-s|\le SL$ throughout the interval of existence.
Hence $U\in[s-1,s+1]\subset[1,\infty)$, and integration of the flux
equation bounds $q$ for this shot. Continuation on both steps gives
existence up to $S$. On every compact set of initial heights these
bounds hold on a common rectangle, so the terminal residual is
continuous there. The integration by parts giving
\eqref{eq:exact-identity} applies to every such shot.

Put $B(x)=g(U(x;s))/g(s)$ and write $\|\cdot\|_\infty$ for the
supremum over $[0,S]$. Let $M=\sup_{[1,\infty)}|\ell|$ and
$H=\sup_{[1,\infty)}|\ell'|$. Integrating the logarithmic derivative
and using $|U-s|\le SL$, we obtain
\begin{align*}
 |\log B(x)|&=\left|\int_s^{U(x;s)}\ell(v)\dd v\right|\le MSL,\\
 |B(x)-1|&\le e^{MSL}-1=O(L),\\
 \left|\frac{g'(U)}{g(s)}-\ell(s)\right|
 &=|\ell(U)B(x)-\ell(s)|\le M|B(x)-1|+HSL=O(L).
\end{align*}
These estimates are uniform for $s\ge2$ and $0\le x\le S$.

Set $k=\Lambda g(s)$, $Z=q/k$ and $A=A_\kappa$. Dividing
\eqref{eq:q-compact} by $k$ gives
\[
 Z_x=A'B-\frac{2c_\kappa L}{\kappa}\mathbf 1_{(1,S)}B.
\]
Since $Z(0)=A(0)=0$ and $\int_0^S|A'|\dd x=2$,
\begin{align*}
 \|Z-A\|_\infty
 &\le \|B-1\|_\infty\int_0^S |A'|\dd x
       +\frac{2c_\kappa L}{\kappa}\int_1^S B\dd x\\
 &\le 2\|B-1\|_\infty+2c_\kappa L\|B\|_\infty=O(L).
\end{align*}
To pass back to the inverse flux with a bound independent of $k$, use
\begin{equation}\label{eq:weighted-flux-error}
 A|\psi(kz)-\psi(kA)|\le4|z-A|
 \qquad(A\ge0,\ z\in\R,\ k>0).
\end{equation}
If $A\le2|z-A|$, the bound follows from $|\psi|\le1$.
Otherwise, every $r$ between $A$ and $z$ satisfies $r>A/2$.
Since $rk|\psi'(kr)|\le1$, the mean-value theorem bounds the
left-hand side by $2|z-A|$.

By \eqref{eq:exact-identity} and \eqref{eq:triangle-integral}, the
error satisfies
\begin{align*}
 \bigl|R_{\Lambda,L,\kappa}(s)-c_\kappa(F_\Lambda(s)-2)\bigr|
 &\le \int_0^S A\left|\frac{g'(U)}{g(s)}-\ell(s)\right|\dd x\\
 &\quad+|\ell(s)|\int_0^S A|\psi(kZ)-\psi(kA)|\dd x
       +\frac{2c_\kappa}{\kappa}\int_1^S|B-1|\dd x .
\end{align*}
The first and last integrals are $O(L)$ by the ratio estimates.
By \eqref{eq:weighted-flux-error}, the middle term is bounded by
$4MS\|Z-A\|_\infty=O(L)$. All constants are independent of $s$
and of $L\in(0,L_*)$, proving \eqref{eq:oscillatory-uniform-error}.
\end{proof}

\begin{theorem}\label{prop:oscillatory}
Let $g(s)=s^2 e^{2s-\cos s}$, fix any $\kappa,\Lambda>0$ and set
$\gamma=2$. For every sufficiently small $L>0$, the weight
\eqref{eq:parameter-family} has length ratio $\kappa$ and strictly
negative mean. For this fixed weight, \eqref{eq:problem} has a
sequence of pairwise distinct positive solutions $u_j$ with
\[
 u_j(0)\longrightarrow+\infty,\qquad
 \min_{0\le t\le T}u_j(t)\longrightarrow+\infty .
\]
\end{theorem}

\begin{proof}
As $s\to\infty$, $g(s)\to\infty$ and
$d(\Lambda g(s))\to1$ by Lemma~\ref{lem:kernel}. Since $\ell$ is
bounded, $F_\Lambda(s)-2=\sin s+o(1)$. Thus for
$s_j=(j+\tfrac12)\pi$ and every sufficiently large integer $j$,
$|F_\Lambda(s_j)-2-(-1)^j|<1/4$.

Choose $0<L<L_*$ with $CL<c_\kappa/4$, using the constants of
Lemma~\ref{lem:tail}. For this single $L$ and every such $j$,
\[
 (-1)^j R_{\Lambda,L,\kappa}(s_j)
 \ge c_\kappa(-1)^j(F_\Lambda(s_j)-2)-CL
 >\tfrac34 c_\kappa-\tfrac14 c_\kappa=\tfrac12 c_\kappa .
\]
The residual is continuous, so each interval $(s_j,s_{j+1})$
contains a zero $\widehat s_j$. Lemma~\ref{lem:tail} and
\eqref{eq:uj-derivative} yield positive Neumann solutions for the
same weight, whose integral is negative by \eqref{eq:mean-computation}.
The heights $u_j(0)=\widehat s_j$ are distinct and tend to infinity.
Since $|u_j'|<1$,
$\min_{[0,T]}u_j\ge\widehat s_j-T\to+\infty$.
\end{proof}

\paragraph*{Funding.}
This research received no specific grant from any funding agency.

\paragraph*{Declaration of competing interest.}
The author declares no competing interests.

\paragraph*{Data availability.}
No external data were used. The numerical results were computed
using the parameters and solver settings given in Example~\ref{ex:numerical}.

\section*{Declaration of generative AI and AI-assisted technologies}

The author used OpenAI Codex to assist with English-language editing
and manuscript organization, and to identify the parameter set used
in Example~\ref{ex:numerical}. The author reviewed and
revised the manuscript, verified the reported shooting residuals and
approximate roots, and takes responsibility for the final article.

\printbibliography[title={References}]
\end{document}